\documentclass{article} \usepackage{amssymb}\usepackage{amscd}
\usepackage[utf8]{inputenc}
\usepackage{amsfonts}
\usepackage{amstext}
\newtheorem{theorem}{Theorem}
\newtheorem{lemma}{Lemma}\newtheorem{claim}{Claim}

\newtheorem{example}{Example}
 
\newtheorem{corollary}{Corollary} \newcommand{\La}{\Lambda}
\newcommand{\R}{{\mathbb R}}  \newcommand{\Z}{{\mathbb Z}} \newcommand{\N}{{\mathbb N}}

\newcommand{\Cc}{{\mathbb C}}

\begin{document}
\title{A Simple Crystalline Measure}\author{Alexander Olevskii  and Alexander Ulanovskii}\date{}\maketitle

\abstract{We prove that every pair of exponential polynomials with imaginary frequencies generates a Poisson-type formula.}

%\begin{large}
\section{Introduction}

Following Yves Meyer  \cite{m16}, %we call an atomic measure $\mu$  a "crystalline measure"             if  the Fourier transform mu^  (in sense of distributions )             is also atomic measure            and the support and the spectrum  of \mu   both are locally finite sets .
we say that   $\mu$ is  a {\it crystalline measure}
            if  $\mu$ is a purely atomic measure supported by a locally finite set
            and its distributional Fourier transform is also a purely atomic measure supported by a locally finite set.
This concept is intimately connected with analysis of adequate mathematical models for the quasicrystal phenomenon.

The classical example of a crystalline measure is given by the Poisson measure$$\mu=\sum_{k\in\Z}\delta_k,$$
            which assigns to each integer  a unit  mass.   This measure, often called the
             Dirac comb, corresponds to crystals.

           The existence of crystalline measures  whose support does not
            contain an infinite arithmetic progression is a non-trivial
            fact.   At present several examples of such measures are known  based
            on different constructions.  The reader can find information on the history of problem  and related results
   updated to 2019 in:  Yves Meyer \cite{m94}, Jeffry Lagarias \cite{l00},  Nir Lev and Alexander Olevskii \cite{lo15}, \cite{lo16} and \cite{lo17},
   Yves Meyer \cite{m16} and \cite{m17}. See also the references in these papers.

               Recently a new progress was obtained by Pavel Kurasov and Peter Sarnak \cite{ks20},
               who  constructed {\it positive }  crystalline measures with
               {\it uniformly discrete} support.

                Even more recently, Yves Meyer has given a geometric approach to construction
              of crystalline measures  based on the concept of {\it curve model set} and
              used inner holomorphic functions, see \cite{m20}. The example from \cite{ks20}
              is obtained in  this work.

             Inspired by  \cite{ks20} and  \cite{m20},  we present an approach to the problem which gives
               simple examples of crystalline measures satisfying the properties mentioned above.

\section{Poisson-type Formulas}
Recall that a set $\La$ of real or complex numbers  is called uniformly discrete (u.d.) if
\begin{equation}
\inf_{\lambda,\lambda'\in\La,\lambda\ne\lambda'}|\lambda-\lambda'|>0.
\end{equation}
%The constant $d(\La)$ is called the separation constant of $\La$.

The classical Poisson summation formula states
$$
\sum_{k\in\Z}{f(k)}=\sum_{k\in\Z}\hat f(k),\quad f\in S(\R),
$$where $S(\R)$ is the Schwartz space and $\hat f$ is the  Fourier transform of $f$.

We present an approach which yields a number of   Poisson-type formulas
\begin{equation}\label{p}
\sum_{\lambda\in\La}c_\lambda{f(\lambda)}=\sum_{s\in S}a_s\hat f(s),
\end{equation}
where $\La$ is a u.d. set in $\Cc$ or $\R$, $S\subset\R$ is a locally finite set satisfying
\begin{equation}\label{sr}
\# S\cap (-r,r)\leq Cr^m,\quad r>1,
\end{equation}with some $C$ and $m$,
and the sequence of coefficients $c_\lambda$ is bounded. If $\La\subset\R$ then the sequence $a_s$ is bounded, too.
The formula holds for certain classes of functions $f$ specified below.

%\begin{definition}

We now introduce three classes of exponential polynomials.
Let $Exp$ denote the set of all exponential polynomials
$$
\psi(z)=\sum_{j=1}^l{b_j}e^{2\pi i \alpha_j z},
$$where   $l\in\N,\alpha_j\in\R$, $b_j\in\Cc$. We denote by $$\La_\psi:=\{z\in\Cc:\psi(z)=0\}$$ the zero set of $\psi.$
One may easily check that for every $\varphi\in Exp$ there is a positive $r$ such that $\La_\varphi$ lies in the strip $|$Im$\,z|<r$.

Let $Exp_s$ denote the subset of $Exp$ of  exponential polynomials $\varphi$ such that   the zeros  of $\varphi$ are simple and the zero set $\La_\varphi$ is  u.d.

Finally, let $Exp_r$ denote  the subset of $Exp_s$ of  exponential polynomials $\varphi$ such that $\La_\varphi\subset\R$.
%\end{definition}
%Denote by $\La_\varphi:=\{z\in \Cc: \varphi(z)=0\}$ the zero set of $\varphi$.

%One may easily check that for every $\varphi\in Exp$ there is a positive $r$ such that $\La_\varphi$ lies in the strip $|$Im$\,z|<r$.

%Denote by $Exp_r$ the subset of $Exp$ of functions $\varphi$ such that $\La_\varphi\subset\R$.

 It turns out that every pair of functions $\psi\in Exp, \varphi\in Exp_s$ satisfying\footnote{Observe that by  a theorem of Ritt \cite{r},  condition $\Lambda_\varphi \subset\La_\psi$ holds if and only if there exists $\varphi_1\in Exp$ such that $\psi(z)=\varphi(z)\varphi_1(z).$}
\begin{equation}\label{subs}\Lambda_\varphi \not\subset\La_\psi.\end{equation}  gives rise to  a Poisson-type formula:

\begin{theorem}
 Given any pair  $\psi\in Exp,\varphi\in Exp_s$ satisfying (\ref{subs}).

(i) Formula  (\ref{p}) holds with $\La=\Lambda_\varphi$,  some set $S\subset\R$ satisfying (\ref{sr}) and a bounded sequence  $c_\lambda$ defined by \begin{equation}\label{coe}c_\lambda=\frac{\psi(\lambda)}{\varphi'(\lambda)},\quad\lambda\in\La.\end{equation}
Formula (\ref{p}) holds for every  entire function $f$ satisfying  $f(x+iy)\in S(\R),$ for every fixed $y$.

(ii) If $\varphi\in Exp_r$, then the sequence $a_s$ is bounded and (\ref{p}) holds for every $f\in S(\R)$.

\end{theorem}

The  set $S$ and coefficients $a_s$ in (\ref{p}) will also be explicitly defined. We also remark that the classes of functions $f$ in Theorem 1 (i) and (ii) can be substantially enlarged, see  the construction below.

We will prove Theorem 1 in sec. 4.

Observe that  Theorem 1 admits an extension: In fact every pair  of exponential polynomials $\psi,\varphi\in Exp$ satisfying (\ref{subs}) generates a Poisson-type formula.
 If   the roots of $\varphi$ are not simple and the set $\La_\varphi$ is not u.d., then the left hand-side of this formula may contain derivatives of $f$ and the coefficients may no longer be bounded. For simplicity of presentation, we will not consider this case.

%\end{remark}

%This can be proved similarly to the proof of Theorem 1 below.

\section{Crystalline Measures}

 %Denote by $Exp_r$ the subset of $Exp_s$ of functions $\varphi$ satisfying $\La_\varphi\subset\R$.

 The following is an immediate corollary of Theorem 1:

 \begin{corollary}\label{c1}
 Given any pair $\varphi\in Exp_r, \psi\in Exp$ satisfying (\ref{subs}). Let $c_\lambda$ be defined in (\ref{coe}). Then
 $$\mu:=\sum_{\lambda\in\La_\varphi}c_\lambda\delta_\lambda$$ is a crystalline measure whose Fourier transform is given by
 $$
 \hat\mu=\sum_{s\in S}a_s\delta_s.
 $$
 \end{corollary}

 Choose  $\psi=\varphi'.$ Then in (\ref{coe}) we have   $c_\lambda=1,\lambda\in\La$. This gives

 \begin{corollary}
 Assume  $\varphi\in Exp_r$. Then
 \begin{equation}\label{cr}\mu:=\sum_{\lambda\in\La_\varphi}\delta_\lambda\end{equation} is a crystalline measure.
 \end{corollary}

We see that  every exponential polynomial $\varphi\in Exp_r$ gives rise to  crystalline measures supported by $\La_\varphi$. A `non-trivial' crystalline measure is a measure which is not  a linear combination  of Poisson measures. Hence, one wishes to find $\varphi\in Exp_r$ such that   $\La_\varphi$ does not contain any arithmetic progression. A simple examples  is given by
% linear combination of two sine or cosine-functions. For simplicity, we consider the function
\begin{equation}\label{var}
\varphi(z):=\sin\pi z+\delta\sin z. %,\quad \varphi_2(z):=\cos\pi z+\delta\cos z.
\end{equation}

\begin{example}
Let $ 0<\delta\leq 1/2$, and $\varphi$ be defined in (\ref{var}). Then

(i) $\La_\varphi=\{k+\delta_k:k\in\Z\}$, where $\delta_k\in [-1/6,1/6]$ satisfy \begin{equation}\label{de}|\sin \pi\delta_k|\leq\delta,\quad k\in\Z;\end{equation}

(ii)   $\La_{\varphi}$   does not contain any arithmetic progression.

%(ii) $\La$ has exactly one zero on each interval $(k-\theta,k+\theta)$, where $\theta=\theta(\delta)\in(0,\pi/2)$ is the angle such that $\sin\theta=\delta.$
\end{example}

We will prove this in sec. 5.

The proof below shows that statements (i) and (ii) in the example hold for every function $$\varphi(z)=\sin\pi z+\sum_{j=1}^nd_j\sin\alpha_j z,\quad n\in\N, d_j\in\R, 0<\alpha_j<\pi,$$
 provided at least one ratio $\pi/\alpha_j$ is irrational and $\delta:=\sum_{j=1}^n|d_j|\leq 1/2.$

 Let $\varphi$ be given in (\ref{var}). By (\ref{de}), we see that  the numbers $\delta_k$ in Example~1  satisfy $\delta_k\to0$ as $\delta\to0$. Therefore, the set $\La_\varphi$ `approaches' the set of integers $\Z$.  By  Corollary 2, this gives

\begin{corollary}
For every $\varepsilon>0$ there is a set  \begin{equation}\label{set}\La=\{k+\delta_k:k\in\Z\},\quad 0\leq |\delta_k|<\varepsilon, k\in\Z,\end{equation} which   does not contain any arithmetic progression, and the corresponding measure in (\ref{cr})  is crystalline.
\end{corollary}

On the other hand, it follows from \cite{f}, that  if
 $\La$ satisfies (\ref{set}) where $\delta_k\to 0,|k|\to\infty,$ and the measure $\mu$ in (\ref{cr})  is crystalline, then $\La=\Z$ and so $\mu$ is the Dirac comb.

\section{Proof of Theorem 1}
Denote by $C$ different positive constants.

Fix two functions
$$
\psi(z)=\sum_{j=1}^l b_je^{2\pi i \beta_jz}\in Exp,\quad \varphi(z)=\sum_{j=1}^m d_je^{2\pi i \alpha_jz}\in Exp_s,
$$
where $l\in \N, m\geq 2$, $\beta_1<...<\beta_l, \alpha_1<...<\alpha_m$. We assume that the coefficients $b_j,d_j$ are different from zero.
We also assume that condition (\ref{subs}) is true.

The proof consists of several steps.

\medskip\noindent
1. Denote by   $\gamma$ the contour which consists of two parallel lines $\gamma_1:=\R-iR_-$ and $\gamma_2:=\R+iR_+$, where the latter line is oriented from $\infty$ to $-\infty$ and the numbers $R_\pm>0$  are defined below.

Denote by $E$ the set of all entire  functions $f$ of {\it finite exponential type} satisfying for every $r>0$ the condition \begin{equation}\label{zz}\sup_{z\in\Cc, |{\rm Im}\,z|<r}|z|^2|f(z)|<\infty.\end{equation}
Clearly, for every $f\in E$  the Fourier transform $\hat f$ is continuous and has compact support.

Consider the integral:
$$
I:=\int_\gamma f(z)\frac{\psi(z)}{\varphi(z)}\,dz=I_1+I_2,\quad I_j:=\int_{\gamma_j}f(z)\frac{\psi(z)}{\varphi(z)}\,dz.
$$

We will use two different methods for calculating $I$.

\medskip\noindent
2. Theorem~3 in \cite{l31} states that when $z$ is uniformly bounded
from the zeros of $\varphi$, then $\varphi(z)$ is uniformly bounded from zero (This can be also  proved somewhat similarly to the proof of Claim 1 below). Since $\La_\varphi$ is u.d., we can find a small $\epsilon>0$ such that every component of the  Minkowski sum  $$\La_\varphi+\{z:|z|\leq\epsilon\}$$ is bounded.
Hence,  we may apply the residue theorem to get
 $$
 I=2\pi i\sum_{\lambda\in\La} \mbox{Res} (f\frac{\psi}{\varphi})(\lambda)=2\pi i\sum_{\lambda\in\La}c_\lambda f(\lambda),
 $$where $c_\lambda$ are defined in (\ref{coe}). Observe that by (\ref{subs}), all coefficients $c_\lambda$ cannot be equal to zero.

 We have to check that $c_\lambda$ are bounded from above. This follows from

 \begin{claim}
 We have $\inf_{\lambda\in\Lambda_\varphi}|\varphi'(\lambda)|>0$.
 \end{claim}

Indeed, set $\varphi_\lambda(z):=\varphi(z+\lambda),\lambda\in\La_\varphi$. Then $\varphi_\lambda(0)=0$.

 Recall that the zeros of $\varphi$ are simple and the zero set $\La_\varphi$ is u.d. Hence, there exists $\epsilon>0$ such that every function
 $\varphi_\lambda$ has exactly one root in the circle $|z|\leq\epsilon.$
     By the   Argument Principle,
    $$
    \oint_{|z|=\epsilon}\frac{\varphi'_\lambda(z)}{\varphi_\lambda(z)}dz=2\pi i,\quad \lambda\in\La_\varphi.
    $$

    Assume the claim is not true: there is a sequence $\lambda_j\in\La_\varphi$ such that $\varphi'(\lambda_j)\to0, j\to\infty$.
    Clearly, there is a subsequence $\lambda_{j(k)}$ such that $\varphi_{\lambda_{j(k)}}$
   converge uniformly on compacts to some exponential polyniomial $\tilde\varphi$,  so that $$
    \oint_{|z|=\epsilon}\frac{\tilde\varphi'(z)}{\tilde\varphi(z)}dz=2\pi i.
    $$On the other hand, clearly $\tilde\varphi$ has a double zero at the origin, and so the above formula contradicts  the   Argument Principle.

\medskip\noindent
3. Let us calculate $I_j, j=1,2.$
Write
$$
\frac{1}{\varphi(z)}=\frac{1}{d_1e^{2\pi i \alpha_1 z}+...+d_me^{2\pi i \alpha_m z}}=\frac{(1/d_1)e^{-2\pi i\alpha_1z}}{1+...+(d_m/d_1)e^{2\pi i (\alpha_m-\alpha_1)z}}.
$$Since $\alpha_j-\alpha_1>0, j>1,$ there exists $R_+>0$ such that we have
$$
\left|(d_2/d_1)e^{2\pi i (\alpha_2-\alpha_1)z}+...+(d_m/d_1)e^{2\pi i (\alpha_m-\alpha_1)z}\right|<1,\quad \mbox{Im}\,z\geq R_+.
$$
Hence, for Im$\,z\geq R_+$ we have
$$
\frac{1}{\varphi(z)}=\frac{e^{-2\pi i\alpha_1z}}{d_1}\sum_{k=0}^\infty(-1)^k\left(\frac{d_2}{d_1}e^{2\pi i (\alpha_2-\alpha_1)z}+...+\frac{d_m}{d_1}e^{2\pi i (\alpha_m-\alpha_1)z}\right)^k.
$$

It easily follows that
there is a discrete set $S_+$ and coefficients $p_s, s\in S_+,$ such that
\begin{equation}\label{1}
\frac{\psi(z)}{\varphi(z)}=\sum_{s\in S_+}p_se^{2\pi i sz}, \quad \mbox{Re}\,z\geq R_+.
\end{equation}
One may  check that $$
S_+\subset \{\beta_1-\alpha_1,...,\beta_l-\alpha_1\}+\bigcup_{k=0}^\infty\bigcup_{n_j\geq0,\sum n_j=k}\{n_2(\alpha_2-\alpha_1)+...+n_m(\alpha_m-\alpha_1)\}.
$$

\begin{claim}
$S_+\subset[\beta_1-\alpha_1,\infty)$ and satisfies (\ref{sr}).
%\begin{equation}\label{ss} \# S_+\cap[\beta_1-\alpha_1,R)\leq C^m l R^m, \quad R\geq 1. \end{equation}
\end{claim}

We omit the simple proof.

Similarly to above, there exists $R_->0$ and a locally finite set $S_-\subset[\alpha_m-\beta_l,\infty)$ and coefficients $q_s$ such that
\begin{equation}\label{2}
\frac{\psi(z)}{\varphi(z)}=\sum_{s\in S_-}q_se^{-2\pi i sz}, \quad \mbox{Re}\,z\leq -R_-.
\end{equation}Again, one may check that $S_-$ satisfies (\ref{sr}).
%$$ \# S_-\cap[\alpha_m-\beta_m,R)\leq C^m l R^m, \quad R\geq 1.$$

Let us now calculate the integrals $I_j$. Using (\ref{1}) and (\ref{2}) we get
$$I_1=\sum_{s\in S_-}q_s\hat f(s),\ I_2=-\sum_{s\in S_+}p_s\hat f(-s), \quad f\in E,
$$where each series  converges absolutely, since it contains only a finite number of elements.

\medskip\noindent 3. Comparing the above calculations of integral $I$, we see that formula (\ref{p}) is valid for every $f\in E$, where
$\La=\La_\varphi$, the sequence
$c_\lambda$ is defined in (\ref{coe}) and bounded, the set  $S:=S_-\cup(-S_+)$ satisfies (\ref{sr}),
$a_s=q_s/2\pi i, s\in S_-$, and $a_s=-p_s/2\pi i, s\in -S_+$.

 Let us prove that if $\La_\varphi \subset\Z$ then the sequence of coefficients $a_s$ is  bounded.

\begin{lemma}
Assume formula (\ref{p}) holds for all functions $f\in E$, where  $\La\subset\R$ is a u.d. set, $S\subset\R$ is a locally finite set and the sequence $c_\lambda$ is  bounded. Then the sequence $a_s$ is bounded, too.
\end{lemma}

Proof. For every  $0<\beta<1$ and $s\in S$ set $$
\psi_{\beta,s}(z):=e^{2\pi i sz}\beta\left(\frac{\sin (\pi\beta z)}{\pi \beta z}\right)^2.
$$It is clear that $\psi_{\beta,s}\in E$ and that the $L^1(\R)$-norm $$L:=\|\psi_{\beta,s}\|_1$$ is finite and {\it does not depend} on $\beta$ and $s$.
The Fourier transform $\hat\psi_{\beta,s}$ is equal to $L$ at the point $s$  and vanishes outside the interval $(s-\beta,s+\beta)$.

Since $S$ is locally finite, we may choose $\beta$ so small that $S\cap(s-\beta,s+\beta)=\{s\}$.
Apply (\ref{p}) with $f=\psi_{\beta,s}$. By Bessel's inequality (see i.e. Proposition 2.7 in \cite{ou1}), the left hand-side of (\ref{p}) admits an estimate
$$
\sum_{\lambda\in\La}|c_\lambda||\psi_{\beta,w}(\lambda)|\leq C,
$$where $C$ does not depend on $\beta$.
The right hand-side of (\ref{p}) contains only one term $a_s \hat\psi_{\beta,s}(s)=La_s$.
Hence, $|a_s|<C, s\in S$, which proves the lemma.

If the set $\La_\varphi$ does not lie on $\R$, a similar argument shows that the coefficients $a_s$ have at most exponential growth.

%: \begin{equation}\label{egr} |a_s|\leq C e^{2\pi Rs},\quad R=\max\{R_-,R_+\}. \end{equation}

\medskip\noindent 4. To prove Theorem 1 (i), it remains to check that (\ref{p}) holds for all entire functions $f$ satisfying $f(x+iy)\in S(\R)$ for every fixed $ y\in\R$. Indeed, the Fourier transform $\hat f $ of such a function satisfies 
$$
|\hat f(t)|\leq C_re^{-r|t|},\quad \mbox{for every } r>0,
$$where $C_r$ depends on $r$. Hence,  both sides of (\ref{p}) absolutely converge and so the proof follows by a usual approximation argument.  In fact, one may check that (\ref{p}) remains true for all  functions $f$ analytic in the strip $|$Im$z|\leq r$ and satisfying (\ref{zz}), where  
 $r=\max\{R_-,R_+\}$.

\medskip\noindent 5. Finally, let $\varphi\in Exp_r$. Then $\La=\La_\varphi\subset\R$.
Theorem 1 (ii) can be deduced from Theorem 1 (i) by a usual approximation argument, taking in mind that  both sides of (\ref{p}) converge absolutely for every $f\in S(\R)$.  Again, one may check that (\ref{p}) remains true for all $f$ satisfying
$$
\sup_{x\in\R}(1+|x|^2)|f(x)|<\infty,\quad \sup_{x\in\R}|x|^{m+2}|\hat f(x)|<\infty.
$$

\section{The Example}
Let $\varphi$ be given in (\ref{var}).

\medskip\noindent
(i) It is clear that $\varphi$ changes the sign  on each interval $(k-\gamma,k+\gamma),k\in\Z,k\ne0,$ where $\gamma=\gamma(\delta)\in(0,1/6]$ is defined by
$$
|\sin\pi\gamma|=\delta.
$$
It remains to  check that it has exactly one zero on every such interval, and has no other zeros.

Denote by $\Gamma_n\subset\Cc$   the square $$\Gamma_n:=\{z=x+iy: \max (|x|,|y|)=n+1/2\},\quad n\in\N.$$ Clearly, $\sin\pi z$ has exactly
$2n+1$ zeros inside $\Gamma_n$. By above,  $\varphi(z)=\sin\pi z+\delta\sin z$ has at least $2n+1$ zeros inside $\Gamma_n$.
It is easy to see that $|\sin\pi z|>\delta|\sin z|$ on $\Gamma_n$, provided $n$ is large enough. So, the result follows from Rouch\'{e}'s theorem.

\medskip\noindent
(ii) Consider the `projection' $\La_1$ of $\La$ onto $[-1/2,1/2)$ defined as
$$
\La_1:=[-1/2,1/2)\cap \left(\bigcup_{k\in\Z}(\La+k)\right).
$$Then $\La_1\subset (-1/6,1/6)$.

Let us check that  $\La$ does not contain any arithmetic progression. Assume, to the contrary, that $\alpha+\beta\Z\subset\La$ for some $\alpha$ and $\beta>0$. If $\beta$ is irrational,   the projection of
$\alpha+\beta\Z$ onto $[-1/2,1/2)$   is  everywhere dense, which is not the case. If $\beta=p/q$ for some $p,q\in\N$, then $\sin\pi z+\delta\sin z$ vanishes on  $\alpha+p\Z$, which is clearly  not possible.

%\end{large}

\medskip

\noindent
A.O.: Tel Aviv University,  P.O. Box 39040, Tel Aviv 6997801, Israel, olevskii@yahoo.com

\medskip

\noindent A.U.: Stavanger University,  4036 Stavanger, Norway, alexander.ulanovskii@uis.no
\end{document}